\documentclass[12pt, notitlepage]{article}
\usepackage{amsmath, amssymb, amsthm}
\usepackage{graphicx}
\newtheorem{prop}{Proposition}[section]

\newtheorem{thm}[prop]{Theorem}
\newtheorem{obs}[prop]{Observation}

\newtheorem{conj}[prop]{Conjecture}
\newtheorem{quest}[prop]{Question}

\newcommand{\SYT}{\operatorname{SYT}}

\author{Joel Brewster Lewis \\
Massachusetts Institute of Technology \\
\texttt{jblewis@math.mit.edu}}
\title{Generating trees and pattern avoidance in alternating permutations}

\begin{document}
\maketitle

\begin{abstract}
We extend earlier work of the same author to enumerate alternating permutations avoiding the permutation pattern $2143$.  We use a generating tree approach to construct a recursive bijection between the set $A_{2n}(2143)$ of alternating permutations of length $2n$ avoiding $2143$ and standard Young tableaux of shape $\langle n, n, n\rangle$ and between the set $A_{2n + 1}(2143)$ of alternating permutations of length $2n + 1$ avoiding $2143$ and shifted standard Young tableaux of shape $\langle n + 2, n + 1, n\rangle$.  We also give a number of conjectures and open questions on pattern avoidance in alternating permutations and generalizations thereof.
\end{abstract}

\section{Introduction}

In \cite{me-big}, the author enumerated the set $A_{2n}(1234)$ of alternating permutations with no four-term increasing subsequence by showing that it is in bijection with the set of standard Young tableaux of shape $\langle n, n, n\rangle$.  One of the two proofs given of this result may be understood as a demonstration that these permutations and tableaux have isomorphic generating trees, that is, that they have the same recursive structure.  (The other proof used a modification of the RSK correspondence first used by Ouchterlony \cite{ouch} in the context of doubly-alternating permutations.)  In this paper, we show that this approach can be extended to enumerate alternating permutations avoiding the pattern $2143$.  In particular, we show that $A_{2n}(2143)$ is in bijection with standard Young tableaux of shape $\langle n, n, n\rangle$ (and so also with $A_{2n}(1234)$) and that $A_{2n + 1}(2143)$ is in bijection with \emph{shifted} standard Young tableaux of shape $\langle n + 2, n + 1, n\rangle$.  We end with a large number of open problems.

In Section \ref{sec:prelims}, we provide some definitions and notation that will be used throughout the paper.  In Section \ref{sec:tableaux}, we describe a generating tree for standard Young tableaux of certain shapes and show that it obeys a simple two-parameter labeling scheme.  In Section \ref{sec:1234}, we briefly review (without proof) one of the main results of \cite{me-big}, namely that alternating permutations of even length with no four-term increasing subsequence share the same generating tree as the tableaux of Section \ref{sec:tableaux}.  In Section \ref{sec:2143}, we show that alternating permutations of even length avoiding the pattern $2143$ also share this generating tree.  In Section \ref{sec:2143odd}, we show that similar methods can be used to enumerate alternating permutations of \emph{odd} length avoiding the pattern $2143$, and that the resulting bijection is with shifted Young tableaux of certain shapes.  In Section \ref{sec:open} we give several conjectures and open problems relating to pattern avoidance in alternating permutations.  In Section \ref{sec:appendix} we provide 
data that forms the basis for several of these conjectures.

\section{Preliminaries}\label{sec:prelims}

In this section, we provide important definitions and notation that will be used throughout the rest of the paper.

\subsection{Permutations, partitions, tableaux}

A \emph{permutation} $w$ of length $n$ is a word containing each of the elements of $[n]:=\{1, 2, \ldots, n\}$ exactly once.  The set of permutations of length $n$ is denoted $S_n$.  Given a word $w = w_1 \cdots w_n$ and a permutation $p = p_1\cdots p_k \in S_k$, we say that $w$ \emph{contains the pattern $p$} if there exists a set of indices $1 \leq i_1 < i_2 < \ldots < i_k \leq n$ such that the subsequence $w_{i_1}w_{i_2}\cdots w_{i_k}$ of $w$ is order-isomorphic to $p$, i.e., such that $w_{i_\ell} < w_{i_m}$ if and only if $p_\ell < p_m$.  Otherwise, $w$ \emph{avoids} $p$.  Given a pattern $p$ and a set $S$ of permutations, we denote by $S(p)$ the set of elements of $S$ that avoid $p$.  For example, $S_n(123)$ is the set of permutations of length $n$ avoiding the pattern $123$, i.e., the set of permutations with no three-term increasing subsequence.  We denote the size of a set $S$ by $|S|$, so $|S_3(123)| = 5$.

A permutation $w = w_1w_2\cdots w_n$ is \emph{alternating} if $w_{1} < w_{2} > w_3 < \ldots$.  (Note that in the terminology of \cite{EC1}, these ``up-down'' permutations are \emph{reverse alternating} while alternating permutations are ``down-up'' permutations.  Luckily, this convention doesn't matter: any pattern result on either set can be translated into a result on the other via \emph{complementation}, i.e., by considering $w^c$ such that $w^c_i = n + 1 - w_i$.  Then results for the pattern $123$ and up-down permutations are replaced by results for $321$ and down-up permutations, and so on.)  We denote by $A_n$ the set of alternating permutations of length $n$.  

A \emph{partition} is a weakly decreasing, finite sequence of nonnegative integers.  We consider two partitions that differ only in the number of trailing zeroes to be the same.  We write partitions in sequence notation, as $\langle \lambda_1, \lambda_2, \ldots, \lambda_n \rangle$.

Given a partition $\lambda = \langle \lambda_1, \lambda_2, \ldots, \lambda_n \rangle$, the \emph{Young diagram} of shape $\lambda$ is the left-justified array of $\lambda_1 + \ldots + \lambda_n$ boxes with $\lambda_1$ in the first row, $\lambda_2$ in the second row, and so on.  We identify each partition with its Young diagram and speak of them interchangeably.  If $\lambda$ is a Young diagram with $m$ boxes, a \emph{standard Young tableau} of \emph{shape} $\lambda$ is a filling of the boxes of $\lambda$ with $[m]$ so that each element appears in exactly one box and entries increase along rows and columns.  We identify boxes in Young diagrams and tableaux using matrix coordinates, so the box in the first row and second column is numbered $(1, 2)$.  We denote by $\SYT(\lambda)$ the set of standard Young tableaux of shape $\lambda$.

\subsection{Generating trees}\label{sec:treeintro}

Given a sequence $\{\Sigma_n\}_{n \geq 1}$ of nonempty sets with $|\Sigma_1| = 1$, a \emph{generating tree} for this sequence is a rooted, labeled tree such that the vertices at level $n$ are the elements of $\Sigma_n$ and the label of each vertex determines the multiset of labels of its children.  In other words, a generating tree is one particular type of recursive structure in which heredity is determined by some local data.  We are particularly interested in generating trees for which the labels are (much) simpler than the objects they are labeling.  In this case, we may easily describe a generating tree by giving the label $L_1$ of the \emph{root vertex} (the element of $\Sigma_1$) and the \emph{succession rule} $L \mapsto S$ that gives the set $S$ of labels of the children in terms of the label $L$ of the parent.  

Beginning with work of Chung, Graham, Hoggatt, and Kleiman \cite{CGHK}, generating trees have been put to good use in the study of pattern-avoiding permutations, notably in the work of West (see, e.g., \cite{west1995, west1996}).  The usual approach has been to consider subtrees of the tree of all permutations given by the rule that $v \in S_{n + 1}$ is a child of $u \in S_n$ exactly when erasing the entry $n + 1$ from $v$ leaves the word $u$.  Because we are interested here in alternating permutations and the permutation that results from erasing the largest entry of an alternating permutation typically is not alternating, this tree is unsatisfactory.  However, the inverse tree, in which we arrive at the children of $u$ by inserting arbitrary values in the last position (rather than inserting the largest value in an arbitrary position), \emph{is} well-suited to the case of alternating permutations.  This motivates the following definitions.

Given a permutation $u \in S_n$ and an element $i \in [n + 1]$, there is a unique permutation $v = v_1 v_2 \cdots v_n v_{n + 1} \in S_{n + 1}$ such that $v_{n + 1} = i$ and the word $v_1 v_2 \cdots v_n$ is order-isomorphic to $u$.  We denote this permutation by $u \leftarrow i$ and refer to it as the \emph{extension} of $u$ by $i$.  In other words, the operation of extending $u$ by $i$ replaces each entry $c \geq i$ in $u$ by $c + 1$ and then attaches $i$ to the end of the result.  For example, $3142 \leftarrow 3 = 41523$.

Given a pattern $p$ and a permutation $w \in S_n(p)$, we say that $c \in [n + 1]$ is \emph{active} or \emph{an active value} for $w$ (with respect to $p$) if $w \rightarrow c$ avoids $p$.\footnote{
This terminology is borrowed from the ``usual'' case, in which a position is said to be active (or is called an ``active site'') if one can insert $n + 1$ into that position while preserving pattern avoidance.}

There is a natural generating tree structure on alternating permutations of even length analogous to the tree on all permutations mentioned above: given an alternating permutation $u$ of length $2n$, its children are precisely the alternating permutations $v$ of length $2n + 2$ such that the prefix of $v$ of length $2n$ is order-isomorphic to $u$.  (Of course, there is also a similar tree for alternating permutations of odd length.)  Since pattern containment is transitive, the set $\bigcup_{n \geq 1} A_{2n}(p)$ of alternating permutations of even length avoiding the pattern $p$ is the set of vertices of a 
subtree.  It is these trees (for the patterns $1234$ and $2143$) that we will consider in Section \ref{sec:trees}.

\section{A generating tree for tableaux}\label{sec:tableaux}

For $n \geq 1$, the collection of standard Young tableaux of shape $\langle n, n, n\rangle$ has a simple associated generating tree: given a tableau $S \in \SYT(n, n, n)$, its children are precisely the tableaux $T \in \SYT(n + 1, n + 1, n + 1)$ such that removing the last column of $T$ leaves a tableau that is order-isomorphic to $S$.  Notice that the shape of the tree below $S$ is determined entirely by the entries of the last column of $S$, and in particular by the entries $S(1, n)$ and $S(2, n)$ (since $S(3, n) = 3n$ for all $S \in \SYT(n, n, n)$).  Thus, we wish to choose a labeling for our tree that captures exactly this information.  Our choice (one of several reasonable options) is to assign to each $S \in \SYT(n, n, n)$ the label $(a, b) = (3n + 1 - S(2, n), 3n + 1 - S(1, n))$.  This provides a two-label generating tree for
\[
\bigcup_{n \geq 1} \SYT(n, n, n),
\]
whose root (the unique standard Young tableau of shape $\langle 1, 1, 1\rangle$) has label $(4 - 2, 4 - 1) = (2, 3)$.  In Proposition 4.7 of \cite{me-big}, the author established (in a more general context, and using an object called ``good tableaux'' rather than the equivalent generating tree we've defined here) the following result.

\begin{prop}\label{prop:rectangulartableaux}
The generating tree that we have just described for 
standard Young tableaux with all columns of length three obeys the rule 
\[
(a, b) \mapsto \{(x, y) \mid 2 \leq x \leq a + 1 \textrm{ and } x + 1 \leq y \leq b + 2\}.
\]
\end{prop}
For example, the root tableau 
$\begin{array}{|c|}
\hline
1  \\
\hline
2  \\
\hline
3 \\
\hline
\end{array} $
has label $(2, 3)$ and has five children,
\[\begin{array}{|c|c|}
\hline
1 & 4\\
\hline
2 & 5\\
\hline
3 & 6 \\
\hline
\end{array}
\, , \; 
\begin{array}{|c|c|}
\hline
1 & 3\\
\hline
2 & 5\\
\hline
4 & 6\\
\hline
\end{array}
\, , \;
\begin{array}{|c|c|}
\hline
1 & 2\\
\hline
3 & 5\\
\hline
4 & 6\\
\hline
\end{array}
\, , \;
\begin{array}{|c|c|}
\hline
1 & 3\\
\hline
2 & 4\\
\hline
5 & 6\\
\hline
\end{array}
\, \textrm{ and }
\begin{array}{|c|c|}
\hline
1 & 2\\
\hline
3 & 4\\
\hline
5 & 6\\
\hline
\end{array}
\, ,
\]
which have labels $(2, 3)$, $(2, 4)$, $(2, 5)$, $(3, 4)$ and $(3, 5)$, respectively.  We restate the proof of this result in our current language. 
\begin{proof}
Choose a standard Young tableau $S \in \SYT(n, n, n)$ with label $(a, b)$ and a child $T \in \SYT(n + 1, n + 1, n + 1)$ of $S$ with label $(x, y)$.  We wish to show that $2 \leq x \leq a + 1$ and $x + 1 \leq y \leq b + 2$.  Expressing these conditions in terms of the entries of $S$ and $T$, we must show
\begin{equation}\label{eq:T1}
S(2, n) + 2 \leq T(2, n + 1) \leq 3n + 2
\end{equation}
and
\begin{equation}\label{eq:T2}
S(1, n) + 1 \leq T(1, n + 1) \leq T(2, n + 1) - 1.
\end{equation}

We first demonstrate that Equation \ref{eq:T2} is valid.  It follows from the relationship between $S$ and $T$ that $S(1, n) \leq T(1, n)$, and since $T$ is a standard Young tableau we have $S(1, n) + 1 \leq T(1, n) + 1 \leq T(1, n + 1) \leq T(2, n + 1) - 1$, as needed.

Now we demonstrate that Equation \ref{eq:T1} is valid.  Since $T$ is a standard Young tableau, we have $T(2, n + 1) < T(3, n + 1) = 3n + 3$.  Thus $T(2, n + 1) \leq 3n + 2$, which is the right half of the desired inequality.  For the left half, we compute that 
\begin{align*}
T(2, n + 1) & = \Big| \{(i, j) \mid T(i, j) < T(2, n + 1)\} \Big| \\
& \geq \Big| \{(i, j) \mid T(i, j) < T(2, n) \textrm{ and } j \leq n\} \cup \{(2, n), (1, n + 1)\} \Big| \\
& = \Big| \{(i, j) \mid S(i, j) < S(2, n)\} \cup \{(2, n), (1, n + 1)\} \Big| \\
& = S(2, n) + 2,
\end{align*}
as desired.

We've now established that Equations \ref{eq:T1} and \ref{eq:T2} hold for all children $T$ of $S$, and so that the labels of the children of $S$ satisfy the claimed inequalities.  Now we must show that $S$ has exactly one child with each of these labels.  

Given a parent tableau $S$ and prescribed values $T(1, n + 1)$ and $T(2, n + 1)$ satisfying Equations \ref{eq:T2} and \ref{eq:T1}, we construct the tableau $T$ in the only way possible: define $T(3, n + 1) = 3n + 3$ and for $1 \leq i \leq 3$ and $1 \leq j \leq n$, set
\[
T(i, j) = 
\begin{cases} 
S(i, j),       & S(i, j) < T(1, n + 1), \\
S(i, j) + 1,   & T(1, n + 1) \leq S(i, j) < T(2, n + 1), \\
S(i, j) + 2,   & T(2, n + 1) \leq S(i, j).
\end{cases}
\]
By construction, the object $T$ that results is of shape $\langle n + 1, n + 1, n + 1\rangle$ and contains each of the elements of $[3n + 3]$ exactly once.  It is not difficult to check that Equations \ref{eq:T1} and \ref{eq:T2} imply that $T$ is increasing along rows and columns and consequently that $T \in \SYT(n + 1, n + 1, n + 1)$, and also that $T$ is a child of $S$.  Thus, $S$ has at least one child with each label satisfying the given conditions.  Finally, it's clear that this $T$ is unique.
\end{proof}

\section{Generating trees for $A_{2n}(1234)$ and $A_{2n}(2143)$}\label{sec:trees}

Recall from Section \ref{sec:treeintro} that for any pattern $p$, there is a natural generating tree for the set $\bigcup_{n} A_{2n}(p)$ of $p$-avoiding alternating permutations of even length: given a $p$-avoiding alternating permutation $u$ of length $2n$, its children are precisely the $p$-avoiding alternating permutations $v$ of length $2n + 2$ such that the prefix of $v$ of length $2n$ is order-isomorphic to $u$.  In this section, we study this tree for the patterns $p = 1234$ and $p = 2143$.

\subsection{$1234$}\label{sec:1234}

One of the main results of \cite{me-big} was that even-length alternating permutations with no four-term increasing subsequence share the generating tree for Young tableaux described in Section \ref{sec:tableaux}.  We now briefly describe this labeling.

Given a $1234$-avoiding alternating permutation $w$, assign to it a label $(a, b)$ where $a$ is the smallest entry in $w$ that is the largest entry in a two-term increasing subsequence of $w$, i.e.,
\[
a = \min \{w_j \mid \exists i < j \textrm{ s.t. } w_i < w_j\},
\]
and $b$ is the number of active values for $w$ with respect to $1234$, i.e., $b$ is the number of choices of $c \in [2n + 1]$ such that $w \leftarrow c$ is $1234$-avoiding.  (Equivalently, we could define $b$ to be the smallest entry in $w$ that the largest term in a three-term increasing subsequence, or $b = 2n + 1$ if $w$ contains no three-term increasing subsequence.)  This labeling has previously been used (e.g., by West \cite{west1995} and Bousquet-M\'elou \cite{kernelgentrees}) in the study of $S_n(1234)$. 

For example, we have for the permutation $w = 27583614 \in A_8(1234)$ that $\{w_j \mid \exists i < j \textrm{ s.t. } w_i < w_j\} = \{3, 4, 5, 6, 7, 8\}$ and that the set of active values for $w$ is $\{1, 2, 3, 4\}$.  Thus $a = 3$ and $b = 4$.

As a second example, we have for the permutation $w = 68372514 \in A_8(1234)$ that $\{w_j \mid \exists i < j \textrm{ s.t. } w_i < w_j\} = \{4, 5, 7, 8\}$ and that $w$ contains no three-term increasing subsequence so all of $[9]$ is active for $w$.  Thus $a = 4$ and $b = 9$.

The root of this tree is the permutation $12 \in A_2(1234)$, with label $(2, 3)$.  Propositions 4.3, 4.4 and 4.6 of \cite{me-big} collectively establish that if $w \in A_{2n}(1234)$ has label $(a, b)$ then the collection of labels of children of $w$ is exactly the set of $(x, y)$ such that $2 \leq x \leq a + 1$ and $x + 1 \leq y \leq b + 2$, each pair occurring with multiplicity one.

We now turn our attention to the pattern $2143$.

\subsection{$2143$}\label{sec:2143}

In this section, we show that alternating permutations of even length avoiding $2143$ have a generating tree isomorphic to those mentioned in preceding sections.

Given a permutation $w \in A_{2n}(2143)$, assign to it a label $(a, b)$ where $a = w_{2n - 1}$ is the next-to-last entry of $w$ and $b$ is the number of active values for $w$ in $[2n + 1]$. 

For example, the permutation $w = 68143527 \in A_8(2143)$ has active values $\{1, 2, 3, 8, 9\}$ (the permutations  $792546381$, $791546382$, $791546283$, $691435278$ and $681435279$ avoid $2143$ while the permutations $691435287$, $791435286$, $791436285$ and $791536284$ contain it).  Thus $a = 2$ and $b = 5$.  

As a second example, the permutation $w = 35462718 \in A_8(2143)$ has active values $\{1, 2, 9\}$.  Thus $a = 1$ and $b = 3$.

We have that the root $12 \in A_2(2143)$ has label $(1, 3)$.  We will show that under this labeling the generating tree for $\bigcup_n A_{2n}(2143)$ obeys the rule
\[
(a, b) \mapsto \{(x, y) \mid 1 \leq x \leq a + 1 \textrm{ and } x + 2 \leq y \leq b + 2\},
\]
and we will use this result to establish an isomorphism between this generating tree and those discussed in the preceding sections.  We break the proof of this result into several smaller pieces: Propositions \ref{prop:onlytheselabels} and \ref{prop:capstone} form the meat of the argument establishing the succession rule, while Propositions \ref{prop:activevalues} and \ref{prop:addfirst} are helpful technical lemmas.  We begin with a simple observation that will be of use in the subsequent proofs.

\begin{obs}\label{usefulobs}
Given a permutation $w \in S_n(2143)$, we have that $c \in [n + 1]$ is not active for $w$ if and only if there exists $i < j < k$ such that $w_j < w_i < c \leq w_k$.
\end{obs}

\begin{prop}\label{prop:activevalues}
If $u \in A_{2n}(2143)$ and $u_{2n - 1} = a$ then $\{1, 2, \ldots, a + 1\}$ are active values for $u$.
\end{prop}
\begin{proof}
Fix $u \in A_{2n}(2143)$ with $u_{2n - 1} = a$, choose $c \leq a + 1$ and let $v = u \leftarrow c \in S_{2n + 1}$.  We wish to show that $v$ avoids $2143$, so suppose otherwise.  Then there exist $i < j < k < 2n + 1$ such that $v_i v_j v_k c$ is an instance of $2143$ in $v$.  We use this (suppositional) instance to construct an instance of $2143$ in $u$; this contradiction establishes that $c$ is active for $u$.  In particular, we show that $u_i u_j u_{2n - 2} u_{2n - 1}$ is an instance of $2143$ in $u$ by showing that $v_i v_j v_{2n - 2} v_{2n - 1}$ is an instance of $2143$ in $v$.  In order to do this, it suffices to show that $j < 2n - 2$ (so that $v_i v_j v_{2n - 2} v_{2n - 1}$ is a subsequence of $v$) and that $v_i < v_{2n - 1}$ (so that this subsequence is order-isomorphic to $2143$).

Since $v_i v_j v_k c$ is an instance of $2143$ and $c \leq a + 1$, we have that $v_j < v_i < c \leq a + 1$ and thus $v_i \leq a = u_{2n - 1} \leq v_{2n - 1}$.  There are at least three entries to the right of $v_i$ in $v$ but only two to the right of $v_{2n - 1}$, so $v_i \neq v_{2n - 1}$ and actually $v_i < v_{2n - 1}$, one of the two conditions we need.  It follows that $v_j < v_i < v_{2n - 1} < v_{2n - 2}$ and similarly $v_j < v_{2n}$, so $v_j$ is smaller than all of $v_{2n - 2}, v_{2n - 1}, v_{2n}$ and $v_{2n + 1}$.  These entries form a suffix of $v$, so $v_j$ must occur at an earlier position in $v$.  That is, we have $j < 2n - 2$, the second necessary condition.  Thus $v_i v_j v_{2n - 2} v_{2n - 1}$ is an instance of $2143$ in $v$ and so $u_i u_j u_{2n - 2} u_{2n - 1}$ is an instance of $2143$ in $u$.  Since $u$ avoids $2143$, this is a contradiction, so actually $v$ avoids $2143$ and $c$ is active for $u$, as desired.
\end{proof}

\begin{prop}\label{prop:onlytheselabels}
If $u \in A_{2n}(2143)$ has label $(a, b)$ and $w$ is a child of $u$ with label $(x, y)$ then $1 \leq x \leq a + 1$ and $x + 2 \leq y \leq b+2$.
\end{prop}
\begin{proof}
Suppose that permutations $u \in A_{2n}$ and $w \in A_{2n + 2}$ have the property that the first $2n$ entries of $w$ are order-isomorphic to $u$, and set $a = u_{2n - 1}$.  If $w_{2n + 1} > a + 1$ then also $w_{2n} > a + 1$ and $w_{2n + 2} > a + 1$, while $w_{2n - 1} = a$.  Thus $w^{-1}(a + 1) \not \in \{2n - 1, 2n, 2n + 1, 2n + 2\}$ and so defining $i = w^{-1}(a + 1)$ we have $i < 2n - 1$.  Then $w_i w_{2n - 1} w_{2n} w_{2n + 1}$ is an instance of $2143$ in $w$.  Taking the contrapositive, if $w$ avoids $2143$ then $w_{2n + 1} \leq a + 1$.  Now fix $u \in A_{2n}(2143)$ with label $(a, b)$ and a child $w \in A_{2n + 2}(2143)$ of $u$ with label $(x, y)$; the preceding argument shows that $x \leq a + 1$.  Clearly $x \geq 1$, so we have proved the first half of our assertion.  We now proceed to bound $y$, the number of active values of $w$.

Define a one-to-one function $f\colon [2n] \to [2n + 2]$ by
\[
f(z) = \begin{cases} 
z, & z < w_{2n + 1} \\
z + 1, & w_{2n + 1} \leq z < w_{2n + 2} \\
z + 2, & w_{2n + 2} \leq z,
\end{cases}
\]
so $f(u_\ell) = w_\ell$ for all $\ell \in [2n]$.  We show that $f$ is a map from nonactive values of $u$ to nonactive values of $w$.  To this end, choose any $c \in [2n]$ that is not active for $u$, and choose $i < j < k$ such that $u_j < u_i < c \leq u_{k}$.  We have $(w_i, w_j, w_k, f(c)) = (f(u_i), f(u_j), f(u_k), f(c))$.  One can easily see that $f$ preserves order, so $w_j < w_i < f(c) \leq w_{k}$ and thus $f(c)$ is not active for $w$.  Thus for each of the $2n - b$ choices of a nonactive value $c$ for $u$ we have a corresponding nonactive value $f(c)$ for $w$ and so $w$ has at most $(2n + 2) - (2n - b) = b + 2$ active values, i.e., $y \leq b + 2$.


Finally, we have by Proposition \ref{prop:activevalues} that $\{1, 2, \ldots, x + 1\}$ are active values for $w$.  We also have that $2n + 3$ is active for $w$ and that $2n + 3 \not\in \{1, 2, \ldots, x + 1\}$, so there are at least $x + 2$ active values for $w$.  Thus $x + 2 \leq y$, which completes the proof of our claim.
\end{proof}

\begin{prop}\label{prop:addfirst}
Suppose $u \in A_{2n}(2143)$ has label $(a, b)$ and active values $s_1 < s_2 < \ldots < s_b$.  If $x \leq a + 1$ and $v = u \leftarrow x$ then $v \in A_{2n + 1}(2143)$ and $v$ has active values $1, 2, \ldots, x, s_{x} + 1, s_{x + 1} + 1, \ldots, s_b + 1$.
\end{prop}

For example, $u = 68143527 \in A_8(2143)$ has label $(2, 5)$ and active values $1$, $2$, $3$, $8$ and $9$.  For $x = 3 \leq 2 + 1$ we have $v = 791546283$ with active values $1$, $2$, $3$, $4 \, (= 3 + 1)$, $9\,(= 8 + 1)$ and $10\,(= 9 + 1)$.

\begin{proof} Choose $u \in A_{2n}(1234)$ with label $(a, b)$ and choose $x \leq a + 1$.  Let $s_1 < \ldots < s_b$ be the active values for $u$ and let $v = u \leftarrow x$.  Since $v_{2n + 1} = x \leq a + 1 \leq v_{2n - 1} + 1$ and $v_{2n - 1} < v_{2n}$, we have $v_{2n + 1} \leq v_{2n}$ and thus actually $v_{2n + 1} < v_{2n}$.  Thus $v$ is alternating.  Also, by Proposition \ref{prop:activevalues}, $x$ is an active value for $u$ and so $v$ avoids $2143$.  This proves the first part of our claim, and we now must show that the set of active values for $v$ is $\{1, 2,  \ldots, x,  s_{x} + 1, s_{x + 1} + 1, \ldots, s_b + 1\}$.  We first show that each of these values is in fact an active value for $v$. 
\begin{description}
\item[Case 1.]  Fix $m \geq x$ so that $s_m \geq x$ is an active value for $u$.  We wish to show that $s_m + 1$ is an active value for $v$.  Let $w = w_1 \cdots w_{2n + 2}$ be the result of extending $v$ by $s_m + 1$, and suppose for sake of contradiction that $w$ contains an instance $w_i w_j w_k w_{\ell}$ of $2143$.  Since $v$ avoids $2143$, we must have $\ell = 2n + 2$, so $w_{\ell} = w_{2n + 2} = s_m + 1$.  As $w_{2n + 2} = s_m + 1 > x = v_{2n + 1}$, extending $v$ by $s_m + 1$ does not change the value of the $(2n + 1)$th entry and thus $w_{2n + 1} = v_{2n + 1} = x$, whence $w_{2n + 2} > w_{2n + 1}$.  Since $w_k > w_{2n + 2}$, we have $k \neq 2n + 1$ and so $w_{2n + 1}$ is not part of our instance of $2143$.  Let $v'$ be the permutation order-isomorphic to $w_1w_2 \cdots w_{2n} w_{2n + 2}$.  Thus $v'_i v'_j v'_k v'_{2n + 1}$ is an instance of $2143$ in $v'$.   However, we also have that $v' = u \leftarrow s_m$ and that $s_m$ is active for $u$.  This is a contradiction, so $w$ cannot contain an instance of $2143$, and we conclude that $s_m + 1$ is an active value for $v$ by definition.

\item[Case 2.]  Fix $c \leq x$ so that $c = s_c$ is an active value for $u$.  We wish to show that $c$ is an active value for $v$.  Let $w = w_1 \cdots w_{2n + 1} w_{2n + 2}$ be the result of extending $v$ by $c$, and suppose for sake of contradiction that $w$ contains an instance $w_i w_j w_k w_{\ell}$ of $2143$.  If $\{k, \ell\} \neq \{2n + 1, 2n + 2\}$ then we can conclude by an argument nearly identical to the previous case, so assume that $k = 2n + 1$ and $\ell = 2n + 2$.  Since $w_i < w_k= w_{2n + 1}$ and $w_j < w_k = w_{2n + 1}$ but $w_{2n} > w_{2n + 1}$, $2n \not \in \{i, j\}$ and so $w_i w_j w_{2n} w_{2n + 2}$ is another instance of $2143$ in $w$.  But now we have an instance not including $w_{2n + 1}$ and so we may proceed as in the previous case.  We conclude that $c$ is an active value for $v$.
\end{description}

Finally, we show that these values are the only active values for $v$.  In particular, we must show that for every $c > x$ that is not active for $u$, $c + 1$ is not active for $v$.  Fix such a $c$.  Since $c$ is not active for $u$, there exist $i < j < k$ such that $u_j < u_i < c \leq u_k$.  We have $v_\ell \leq u_\ell + 1$ for all $\ell \in [2n]$, so $v_j < v_i < c + 1$.  Since $u_k > c - 1 \geq x$, we have $v_k = u_k + 1$.  Thus $v_j < v_j < c + 1 \leq u_k + 1 = v_k$ and so $c + 1$ is not an active value for $v$.  We conclude that the active values for $v$ are exactly $1, 2, \ldots, x, s_{x} + 1, s_{x + 1} + 1, \ldots, s_b + 1$, as claimed.  
\end{proof}

\begin{prop}\label{prop:capstone}
If $u \in A_{2n}(2143)$ has label $(a, b)$ and $1 \leq x \leq a + 1$, $x + 2 \leq y \leq b + 2$, then there is a unique child $w \in A_{2n + 2}(2143)$ of $u$ with label $(x, y)$.
\end{prop}
\begin{proof}
Choose a permutation $u \in A_{2n}(2143)$ with label $(a, b)$ and choose $(x, y)$ such that $1 \leq x \leq a + 1$ and $x + 2 \leq y \leq b + 2$.  We will construct a child of $u$ with label $(x, y)$. 

Let $s_1 < \ldots < s_b$ be the active values for $u$.  Define $v = v_1 v_2 \cdots v_{2n} v_{2n + 1}$ and $w = w_1 w_2 \cdots w_{2n} w_{2n + 1} w_{2n + 2}$ by $v = u \leftarrow x$ and $w = v \leftarrow (s_{b + 2 + x - y} + 1)$.  We claim that $w$ is the desired permutation.  We must show that $w$ belongs to $A_{2n + 2}(2143)$ and that its label really is $(x, y)$. 

It follows from Proposition \ref{prop:addfirst} that $v \in A_{2n + 1}(2143)$.  Because $y \leq b + 2$, we have $b + 2 + x - y \geq x$ and so $s_{b + 2 + x - y} + 1 > s_{x} = x = v_{2n + 1}$.  Thus $w$ is alternating.  We also have from Proposition \ref{prop:addfirst} that $s_{b + 2 + x - y} + 1$ is an active value for $v$, so $w \in A_{2n + 2}(2143)$.  We have left to show that $w$ has label $(x, y)$. 

Since $w_{2n + 2} > v_{2n + 1}$, extending $v$ by $w_{2n + 2}$ leaves the value of the entry in the $(2n + 1)$th position unchanged and so $w_{2n + 1} = v_{2n + 1} = x$.  It remains to show that $w$ has exactly $y$ active values.  We claim that the active values for $w$ are precisely $1 < 2 < \ldots < x + 1 < s_{b + 2 + x - y} + 2 < s_{b + 3 + x - y} + 2 < \ldots < s_b + 2$.  By Proposition \ref{prop:activevalues}, we already know that $1, 2, \ldots, x + 1$ are active for $w$.  

For any fixed $i$ such that $b + 2 + x - y \leq i \leq b$, let $z = w \leftarrow (s_i + 2)$.  We wish to show that $z$ avoids $2143$.  We proceed by the same argument as in Case 1 of the proof of Proposition \ref{prop:addfirst}: since $s_i + 2 > w_{2n + 2}$, we have $z_{2n + 2} = w_{2n + 2} < z_{2n + 3}$.  It follows that $z_{2n + 2}$ and $z_{2n + 3}$ cannot be part of the same $2143$ pattern in $z$.  Thus $z$ contains an instance of $2143$ if and only if $z' = z_1 \cdots z_{2n + 1} z_{2n + 3}$ does.  However, $z'$ is order-isomorphic to $v \leftarrow (s_i + 1)$.  Since $s_i + 1$ is an active value for $v$, $z'$ avoids $2143$ and so $z$ avoids $2143$.  Thus $s_i + 2$ is an active value for $w$.  There are $b - (b + 2 + x - y) + 1 = y - x - 1$ such active values.  We must show that there are no other active values for $w$ larger than $x + 1$.

For any fixed $c > s_{b + 2 + x - y} + 2$ that is not of the form $s_i + 2$, we wish to show that $w \leftarrow c$ contains $2143$.  We proceed by the same arguments that follow Case 2 of the proof of Proposition \ref{prop:addfirst}: if $v_iv_jv_k$ are entries of $v$ that can be used to form a $2143$ pattern when $v$ is extended by $c - 1$ then $w_iw_jw_k$ can be used to form a $2143$ pattern when $w$ is extended by $c$.

Finally, for any fixed $c$ such that $x + 1 < c < s_{b + 2 + x - y} + 2$, let $z = w \leftarrow c$.  We wish to show that $z$ contains $2143$.  We have that $z_{2n + 1} = x$, $z_{2n + 2} = s_{b + 2 + x - y} + 2 > c$, and $z_{2n + 3} = c > x + 1$, and we know that there exists $i < 2n + 1$ such that $z_i = x + 1$.  Then $z_i z_{2n + 1}z_{2n + 2} z_{2n + 3} = (x + 1)x(s_{b + 2 + x - y} + 2)c$ is an instance of $2143$ in $z$, so $c$ is not active for $w$.

The preceding four paragraphs account for all elements $[2n + 3]$.  We've shown that exactly $(x + 1) + (y - x - 1) = y$ of these values are active for $w$.  Putting everything together, we have that $w \in A_{2n + 2}(2143)$ is a child of $u$ with $y$ active values and $w_{2n + 1} = x$, i.e., $w$ has label $(x, y)$.  In fact, it follows from our proof that every child of $u$ has a distinct label: we've exhausted the possible pairs of values for $w_{2n + 1}$, $w_{2n + 2}$ such that $w \in A_{2n + 2}(2143)$ is a child of $u$.
\end{proof}

\begin{thm}\label{thm:evenformula}
For all $n \geq 1$ we have
\[
|A_{2n}(2143)| = |\SYT(n, n, n)| = |A_{2n}(1234)| = \frac{2 \cdot (3n)!}{n!(n + 1)!(n + 2)!}.
\]
\end{thm}

\begin{proof}
Proposition \ref{prop:rectangulartableaux} shows that the generating tree for $\bigcup_{n \geq 1}\SYT(n, n, n)$ has root $(2, 3)$ and rule
\[
(a, b) \mapsto \{(x, y) \mid 2 \leq x \leq a + 1 \textrm{ and } x + 1 \leq y \leq b + 2\}
\]
while Propositions \ref{prop:onlytheselabels} and \ref{prop:capstone} together show that the generating tree for $\bigcup_{n \geq 1}A_{2n}(2143)$ has root $(1, 3)$ and rule
\[
(a, b) \mapsto \{(x, y) \mid 1 \leq x \leq a + 1 \textrm{ and } x + 2 \leq y \leq b + 2\}.
\]
These two trees are isomorphic: replacing each label $(a, b)$ in the first tree with $(a - 1, b)$ results in the second tree.  Thus, there is a recursive bijection between $A_{2n}(2143)$ and $\SYT(n, n, n)$, and we have the first claimed equality.  The work in \cite{me-big} summarized in Section \ref{sec:1234} gives the second equality.  Finally, applying the hook-length formula (see, e.g., \cite[Chapter 7.21]{EC2}) gives the third equality.
\end{proof}

\section{Generating tree for $A_{2n + 1}(2143)$}\label{sec:2143odd}

Turning out attention to alternating permutations of odd length, we find that the generating trees of $A_{2n + 1}(1234)$ and $A_{2n + 1}(2143)$ are not isomorphic.  Indeed, the two sequences enumerate differently: all sixteen alternating permutations of length five avoid $1234$, but only twelve of them avoid $2143$.  Although this initially seems disappointing, it turns out that we can still use the methods of the preceding section to enumerate $2143$-avoiding alternating permutations of odd length.

As in \cite{me-big}, it is convenient to consider the set $A'_{2n + 1}$ of \emph{down-up} alternating permutations of odd length rather than up-down alternating permutations; note that results in either case may be translated into results in the other via reverse-complementation.  Arguments very similar to those of Section \ref{sec:2143} show that if we associate to the permutation $w \in A'_{2n + 1}(2143)$ the label $(a, b)$, where $a = w_{2n}$ and $b$ is the number of active values for $w$ then the generating tree for $\bigcup_{n \geq 0} A'_{2n + 1}(2143)$ is as follows:
\begin{description}
\item[Root:] $(0, 2)$
\item[Rule:] $(a, b) \mapsto \{(x, y) \mid 1 \leq x \leq a + 1 \textrm{ and } x + 2 \leq y \leq b + 2\}$
\end{description}
Here the root permutation $1 \in A'_1(2143)$ has label $(0, 2)$.  Its two children $213, 312 \in A'_3(2143)$ have labels $(1, 3)$ and $(1, 4)$, respectively, and so on.  

We seek to enumerate $A'_{2n + 1}(2143)$ by aping our approach for even-length permutations, i.e., by finding a family of objects with isomorphic generating tree that we already know how to enumerate.  In the case at hand, these objects turn out to be \emph{shifted standard Young tableaux} (henceforward SHSYT) of shape $\langle n + 2, n + 1, n \rangle$.  (See \cite{shiftedhook} or \cite[Chapter 10]{skewbook} for definitions, etc.)

Given a SHSYT $T$ of shape $\langle n + 2, n + 1, n\rangle$, assign to it the label $(a, b)$, where $a = 3n + 4 - T(n + 2, 2)$ and $b = 3n + 4 - T(n + 2, 1)$.  (Note that for $n \geq 1$, we have $T(n + 2, 3) = 3n + 3$ for every SHSYT of this shape, so this label captures all the information we need to reconstruct the last column of $T$.)  Then the root of the tree is the unique SHSYT \, \parbox{.5in}{\begin{picture}(0,0)%
\includegraphics{shtab1.pstex}%
\end{picture}%
\setlength{\unitlength}{3947sp}%
\begingroup\makeatletter\ifx\SetFigFont\undefined%
\gdef\SetFigFont#1#2#3#4#5{%
  \reset@font\fontsize{#1}{#2pt}%
  \fontfamily{#3}\fontseries{#4}\fontshape{#5}%
  \selectfont}%
\fi\endgroup%
\begin{picture}(473,474)(964,-1648)
\put(1051,-1336){\makebox(0,0)[lb]{\smash{{\SetFigFont{12}{14.4}{\familydefault}{\mddefault}{\updefault}{$1$}%
}}}}
\put(1275,-1336){\makebox(0,0)[lb]{\smash{{\SetFigFont{12}{14.4}{\familydefault}{\mddefault}{\updefault}{$2$}%
}}}}
\put(1275,-1561){\makebox(0,0)[lb]{\smash{{\SetFigFont{12}{14.4}{\familydefault}{\mddefault}{\updefault}{$3$}%
}}}}
\end{picture}%
}
of shape $\langle 2, 1\rangle$, which has label $(1, 2)$.  Its children are the two SHSYT 
\begin{center}
\parbox{.75in}{\begin{picture}(0,0)%
\includegraphics{shtab2.pstex}%
\end{picture}%
\setlength{\unitlength}{3947sp}%
\begingroup\makeatletter\ifx\SetFigFont\undefined%
\gdef\SetFigFont#1#2#3#4#5{%
  \reset@font\fontsize{#1}{#2pt}%
  \fontfamily{#3}\fontseries{#4}\fontshape{#5}%
  \selectfont}%
\fi\endgroup%
\begin{picture}(699,699)(4189,-2248)
\put(4726,-2161){\makebox(0,0)[lb]{\smash{{\SetFigFont{12}{14.4}{\familydefault}{\mddefault}{\updefault}{$6$}%
}}}}
\put(4726,-1936){\makebox(0,0)[lb]{\smash{{\SetFigFont{12}{14.4}{\familydefault}{\mddefault}{\updefault}{$5$}%
}}}}
\put(4726,-1711){\makebox(0,0)[lb]{\smash{{\SetFigFont{12}{14.4}{\familydefault}{\mddefault}{\updefault}{$4$}%
}}}}
\put(4501,-1711){\makebox(0,0)[lb]{\smash{{\SetFigFont{12}{14.4}{\familydefault}{\mddefault}{\updefault}{$2$}%
}}}}
\put(4501,-1936){\makebox(0,0)[lb]{\smash{{\SetFigFont{12}{14.4}{\familydefault}{\mddefault}{\updefault}{$3$}%
}}}}
\put(4276,-1711){\makebox(0,0)[lb]{\smash{{\SetFigFont{12}{14.4}{\familydefault}{\mddefault}{\updefault}{$1$}%
}}}}
\end{picture}%
} and \hspace{.15in} \parbox{.6in}{\begin{picture}(0,0)%
\includegraphics{shtab3.pstex}%
\end{picture}%
\setlength{\unitlength}{3947sp}%
\begingroup\makeatletter\ifx\SetFigFont\undefined%
\gdef\SetFigFont#1#2#3#4#5{%
  \reset@font\fontsize{#1}{#2pt}%
  \fontfamily{#3}\fontseries{#4}\fontshape{#5}%
  \selectfont}%
\fi\endgroup%
\begin{picture}(699,698)(3814,-1873)
\put(4351,-1786){\makebox(0,0)[lb]{\smash{{\SetFigFont{12}{14.4}{\familydefault}{\mddefault}{\updefault}{$6$}%
}}}}
\put(4351,-1562){\makebox(0,0)[lb]{\smash{{\SetFigFont{12}{14.4}{\familydefault}{\mddefault}{\updefault}{$5$}%
}}}}
\put(3901,-1337){\makebox(0,0)[lb]{\smash{{\SetFigFont{12}{14.4}{\familydefault}{\mddefault}{\updefault}{$1$}%
}}}}
\put(4126,-1337){\makebox(0,0)[lb]{\smash{{\SetFigFont{12}{14.4}{\familydefault}{\mddefault}{\updefault}{$2$}%
}}}}
\put(4351,-1337){\makebox(0,0)[lb]{\smash{{\SetFigFont{12}{14.4}{\familydefault}{\mddefault}{\updefault}{$3$}%
}}}}
\put(4126,-1562){\makebox(0,0)[lb]{\smash{{\SetFigFont{12}{14.4}{\familydefault}{\mddefault}{\updefault}{$4$}%
}}}}
\end{picture}%
}
\end{center}
of shape $\langle 3, 2, 1\rangle$, which have labels $(2, 3)$ and $(2, 4)$, respectively.  In subsequent layers of the tree, the succession rule is identical to the rule for rectangular tableaux; indeed, the subtree below an SHSYT of shape $\langle n + 2, n + 1, n\rangle$ depends only on its last column (or equivalently, on its associated label), and so might as well be the subtree of an SYT of shape $\langle (n + 1)^3\rangle$ with the same label.

It follows immediately that $|A'_{2n + 1}(2143)|$ (and so also $|A_{2n + 1}(2143)|$) is the number of SHSYT of shape $\langle n + 2, n + 1, n\rangle$.  As in the case of standard Young tableaux, there is a simple hook-length formula for SHSYT (see e.g. \ \cite{shiftedhook} or \cite[pp. 187-190]{skewbook}).  In our case, it gives the following result.
\begin{prop}\label{prop:oddformula}
For $n \geq 0$ we have
\[
|A_{2n + 1}(2143)| = \frac{2 (3n + 3)!}{n!(n + 1)!(n + 2)! (2n + 1)(2n + 2)(2n + 3)}.
\]
\end{prop}

\section{Open problems}\label{sec:open}

In this section we pose a number of open enumerative problems related to pattern avoidance in alternating permutations.

\subsection{Other equivalences for patterns of length 4}

If permutations $p$ and $q$ satisfy $|A_{2n}(p)| = |A_{2n}(q)|$ for all $n \geq 1$, we say that $p$ and $q$ are \emph{equivalent} for even-length alternating permutations.  Note that if $p = p_1 \cdots p_k$ and $q = (k + 1 - p_k)(k + 1 - p_{k - 1}) \cdots (k + 1 - p_1)$ (i.e., $p$ and $q$ are reverse-complements) then $p$ and $q$ are equivalent for even-length alternating permutations: for every $n$, the operation of reverse-complementation is a bijection between $A_{2n}(p)$ and $A_{2n}(q)$.  Pairs of patterns that are equivalent for this reason are said to be \emph{trivially equivalent}.  Similarly, if $|A_{2n+1}(p)| = |A_{2n+1}(q)|$ for all $n \geq 0$, we say that $p$ and $q$ are equivalent for odd-length alternating permutations, and if $p$ is the reverse of $q$ then we say they are trivially equivalent.

Numerical data (see Section \ref{sec:appendix}) suggest the following conjecture.

\begin{conj}\label{conj:evenclass}
We have $|A_{2n}(p)| = |A_{2n}(1234)| \; (= |A_{2n}(2143)|)$ for all $n \geq 1$ and every $p \in \{1243, 2134, 1432, 3214, 2341, 4123, 3421, 4312\}$.
\end{conj}
Observe that these eight patterns come in four pairs of trivially equivalent patterns.  The results of West \cite{west1995} and computer investigations of short permutations suggest that some of these equivalences may be susceptible to a generating-tree attack.  In particular, the generating trees for alternating permutations avoiding $1243$ or $2134$ may be isomorphic to the generating tree discussed in Sections \ref{sec:tableaux} and \ref{sec:trees}.

For alternating permutations of even length, the only other possible equivalences not ruled out by numerical data are captured by the following conjecture.
\begin{conj}
We have $|A_{2n}(3142)| = |A_{2n}(3241)| = |A_{2n}(4132)|$ and $|A_{2n}(2413)| = |A_{2n}(1423)| = |A_{2n}(2314)|$ for all $n \geq 1$.
\end{conj}
In both cases, the second of the two equalities is a trivial equivalence.  

For odd-length alternating permutations, computational data suggest the following conjectures.
\begin{conj}\label{conj:1234oddclass}
We have $|A_{2n + 1}(p)| = |A_{2n + 1}(1234)|$ for all $n \geq 0$ and every $p \in \{2134, 4312, 3214, 4123\}$.
\end{conj}
We also have the trivial equivalence $|A_{2n + 1}(1234)| = |A_{2n + 1}(4321)|$.  The equivalence between $4321$ and $4312$ may be amenable to generating tree methods.
\begin{conj}\label{conj:2143oddclass}
We have $|A_{2n + 1}(p)| = |A_{2n + 1}(2143)|$ for all $n \geq 0$ and every $p \in \{1243, 3421, 1432, 2341\}$.
\end{conj}
We also have the trivial equivalence $|A_{2n + 1}(2143)| = |A_{2n + 1}(3412)|$.  The equivalence between $3412$ and $3421$ may be amenable to generating tree methods.

The only other possible equivalence for odd-length alternating permutations not ruled out by data is captured by the following conjecture.

\begin{conj}
The permutations $2314, 4132, 2413, 3142, 1423$ and $3241$ are equivalent for odd-length alternating permutations.
\end{conj}

Other than the cases covered by Conjectures \ref{conj:evenclass}, \ref{conj:1234oddclass} and \ref{conj:2143oddclass}, none of the sequences $\{|A_{2n}(p)|\}_n$ or $\{|A_{2n + 1}(p)|\}_n$ for $p \in S_4$ are recognizable to the present author (and in particular they do not appear in the OEIS \cite{OEIS}).  Numerical data (see Section \ref{sec:appendix}) rule out simple product formulas similar to those of Theorem \ref{thm:evenformula} and Proposition \ref{prop:oddformula}.  However, in some cases it may be possible to give a generating tree and so perhaps to adapt the method of \cite{kernelgentrees} to find generating functions or even closed formulas.

\subsection{Other problems}

The conjectures of the preceding section concern two equivalence relations on patterns of length four: equivalence for even-length alternating permutations and equivalence for odd-length alternating permutations.  It happens that these relations are both refinements of the usual (Wilf-)equivalence relation for all permutations (see, for example, \cite{StankovaWest}).  (Note that this is true on account of numerical data, regardless of the truth of any of the preceding conjectures.)  This suggests the following conjecture.

\begin{conj}\label{conj:subequivalence}
If permutations $p$ and $q$ are equivalent for alternating permutations of either parity then $p$ and $q$ are Wilf-equivalent for all permutations.
\end{conj}

More broadly, we can ask the following question.

\begin{quest}
Are there any large families of patterns that can be shown to be equivalent for alternating permutations (of either parity)?
\end{quest}

The work in \cite{me-big} suggests two possible generalizations of alternating permutations for purposes of pattern avoidance.  One of these (denoted $\mathcal{L}_{n, k}$ in \cite{me-big}; it consists of the reading words of standard Young tableaux of the ``thickened staircase'' shape $\langle n + k - 1, n + k - 2, \ldots, k + 1, k \rangle / \langle n - 1, n - 2, \ldots, 1\rangle$) does not seem to give rise to nontrivial equivalences.  However, an alternative (and perhaps more natural) generalization is the set $\operatorname{Des}_{n, k}$ of permutations of length $n$ with descent set $\{k, 2k, \ldots\}$.  

\begin{quest}
Is $\operatorname{Des}_{n, k}$ a ``good'' context to study pattern avoidance?  In particular, are there any pairs or families of patterns that can be shown to be Wilf-equivalent for these permutations?  Is Conjecture \ref{conj:subequivalence} valid in this context?
\end{quest}
In analogy with the case of alternating permutations, it is natural to consider separate cases depending on the congruence class of $n$ modulo $k$.  Finally, we note that it may also be fruitful to consider permutations whose descent set is contained in (rather than ``is equal to'') $\{k, 2k, \ldots\}$, a case on which our work here does not touch.

\newpage

\section{Appendix}\label{sec:appendix}

The tables of data that follow form the basis for several of the conjectures in the preceding section.  They were generated by brute-force computer enumeration.  

Tables \ref{table:oddaltS4} and \ref{table:evenaltS4} give the number of alternating permutations avoiding patterns of length four, grouped by conjectural equivalence.

\bigskip

\begin{table}[h]
\begin{tabular}{|l||cccccc|}
\hline
Patterns	& 1 & 3 & 5 & 7 & 9 & 11 \\
\hline\hline
$(1234, 4321)$, 
$(2134, 4312)$, 
$(3214, 4123)$
	        & 1 & 2 & 16 & 168 & 2112 & 30030 \\
\hline
$(2143, 3412)$,
$(1243, 3421)$, 
$(1432, 2341)$
		& 1 & 2 & 12 & 110 & 1274 & 17136 \\
\hline
$(2314, 4132)$,
$(2413, 3142)$,
$(1423, 3241)$
		& 1 & 2 & 12 & 106 & 1138 & 13734 \\
\hline
$(1324, 4231)$
                & 1 & 2 & 12 & 110 & 1285 & 17653 \\
\hline
$(1342, 2431)$
                & 1 & 2 & 12 & 108 & 1202 & 15234 \\
\hline
$(3124, 4213)$
                & 1 & 2 & 16 & 168 & 2072 & 28298 \\
\hline
\end{tabular}
\caption{Values of $A_{n}(p)$ for odd $n$ and $p \in S_4$.  Parentheses indicate trivial equivalences.}
\label{table:oddaltS4}
\end{table}

\begin{table}[h]
\begin{tabular}{|p{2.25in}||cccccc|}
\hline
Patterns	& 2 & 4 & 6 & 8 & 10 & 12 \\
\hline\hline
$1234$,
$(1243, 2134)$, 
$(1432, 3214)$, 
$2143$,
$(2341, 4123)$, 
$(3421, 4312)$
					& 1 & 5 & 42 & 462 & 6006 & 87516 \\
\hline
$3142$,
$(3241, 4132)$
					& 1 & 5 & 42 & 444 & 5337 & 69657\\
\hline
$(1423, 2314)$,
$2413$
					& 1 & 4 & 28 & 260 & 2844 & 34564 \\
\hline
$3412$
		& 1 & 4 & 29 & 290 & 3532 & 49100 \\
\hline
$1324$
		& 1 & 4 & 29 & 292 & 3620 & 51866 \\
\hline
$(1342, 3124)$
		& 1 & 5 & 42 & 453 & 5651 & 77498 \\
\hline
$(2431, 4213)$
		& 1 & 5 & 42 & 454 & 5680 & 78129 \\
\hline
$4231$
		& 1 & 5 & 42 & 462 & 6070 & 90686 \\
\hline
$4321$
		& 1 & 5 & 61 & 744 & 10329 & 157586 \\
\hline
\end{tabular}
\caption{Values of $A_{n}(p)$ for even $n$ and $p \in S_4$.  Parentheses indicate trivial equivalences.}
\label{table:evenaltS4}
\end{table}

 \newpage

Tables \ref{table:oddaltS5} and \ref{table:evenaltS5} give the number of alternating permutations avoiding certain patterns of length five.  Only those patterns that might potentially have nontrivial equivalences are included, and they are grouped by these potential equivalences.

\bigskip

\begin{table}[h]
\begin{center}
\begin{tabular}{|p{2.25in}||cccccc|}
\hline
Patterns	& 1 & 3 & 5 & 7 & 9 & 11 \\
\hline\hline
$(12534, 43521)$,
$(21534, 43512)$
                & 1 & 2 & 16 & 243 & 5291 & 144430 \\
\hline
$(12453, 35421)$,
$(21453, 35412)$
                & 1 & 2 & 16 & 243 & 5307 & 146013 \\
\hline
(12354, 45321),
(12543, 34521),
(15432, 23451),
(21354, 45312),
(21543, 34512),
(32154, 45123)
                & 1 & 2 & 16 & 243 & 5330 & 148575 \\
\hline
$(12435, 53421)$,
$(21435, 53412)$
                & 1 & 2 & 16 & 243 & 5330 & 148764 \\
\hline
(12345, 54321),
(21345, 54312),
(32145, 54123),
(43215, 51234)
                & 1 & 2 & 16 & 272 & 6531 & 194062 \\
\hline
\end{tabular}
\end{center}
\caption{Selected values of $A_n(p)$ for odd $n$ and $p \in S_5$.  Parentheses indicate trivial equivalences.  All possible nontrivial equivalences are among the permutations in this table.}
\label{table:oddaltS5}
\end{table}

\begin{table}[h!]
\begin{center}
\begin{tabular}{|p{2.25in}||cccccc|}
\hline
Patterns	& 2 & 4 & 6 & 8 & 10 & 12 \\
\hline\hline
$(12534, 23145)$,
$(21534, 23154)$
                & 1 & 5 & 56 & 997 & 23653 & 679810 \\
\hline
$(34512, 45123)$,
$45312$
                & 1 & 5 & 56 & 1004 & 24310 & 724379 \\
\hline
$(12435, 13245)$,
$(13254, 21435)$
                & 1 & 5 & 56 & 1004 & 24336 & 727807 \\
\hline
$(12453, 31245)$,
$(21453, 31254)$
                & 1 & 5 & 61 & 1194 & 30802 & 953088 \\
\hline
$\phantom{(}12345$,
$21354$,\phantom{)}
$(12354, 21345)$, 
$(12543, 32145)$,
$(15432, 43215)$,
$(21543, 32154)$,
$(23451, 51234)$,
$(34521, 54123)$, \,
$(45321, 54312)$
                 & 1 & 5 & 61 & 1194 & 30945 & 970717 \\
\hline
\end{tabular}
\end{center}
\caption{Selected values of $A_n(p)$ for even $n$ and $p \in S_5$.  Parentheses indicate trivial equivalences.  All possible nontrivial equivalences are among the permutations in this table.}
\label{table:evenaltS5}
\end{table}

 \newpage

Table \ref{table:des3S4S5} gives the number of permutations of length $3n$ with descent set $\{3, 6, \ldots\}$ that avoid certain patterns of length four and five.  
Only those patterns that might potentially have nontrivial equivalences are included.

\bigskip

\begin{table}[h]
\begin{tabular}{|l||cccc|}
\hline
Patterns	& 3 & 6 & 9 & 12 \\
\hline\hline
$2413$,
$(1423, 2314)$
                & 1 & 9 & 153 & 3465 \\
\hline
$(1243, 2134)$,
$(2341, 4123)$
                & 1 & 9 & 153 & 3579 \\
\hline
$3142$,
$(3241, 4132)$
                & 1 & 19 & 642 & 27453 \\
\hline
$2143$,
$4231$,
$(1432, 3214)$,
$(3421, 4312)$
                & 1 & 19 & 642 & 29777 \\
\hline\hline
(12354, 21345),
(23451, 51234)
                & 1 & 19 & 887 & 66816 \\
\hline
(15243, 32415),
(35241, 52413)
                & 1 & 19 & 1077 & 102051 \\
\hline
(12543, 32154),
(34521, 54123)
                & 1 & 19 & 1134 & 114621 \\
\hline
$21354$,
$52341$
                & 1 & 19 & 1134 & 115515 \\
\hline
\parbox[t]{3.48in}{
(15432, 43215),
(21543, 32154),
(25431, 53214),
(31542, 42153),
(32541, 52143),
(35421, 54213),
(41532, 43152),
(42531, 53142),
(43251, 51432),
(43521, 54132),
(45321, 54312),
(52431, 53241),
(53421, 54231)}
                & 1 & 19 & 1513 & 211425 \\
\hline
\end{tabular}
\caption{Selected values of $\operatorname{Des}_{n,3}(p)$ for $n$ divisible by $3$ and $p \in S_4$ or $S_5$.  Parentheses indicate trivial equivalences.  All possible nontrivial equivalences are among the permutations in this table.}
\label{table:des3S4S5}
\end{table}

\bibliographystyle{plain}
\bibliography{2143bib}{}

\end{document}